\documentclass[a4paper,12pt,reqno]{amsart}

\usepackage{amsmath,amssymb,amsthm}
\usepackage{latexsym}
\usepackage{color}
\usepackage{graphicx}
\usepackage{mathrsfs}
\usepackage{enumerate}
\usepackage[abbrev]{amsrefs}
\usepackage[T1]{fontenc}
\usepackage{mathtools}
\mathtoolsset{showonlyrefs=true}

\allowdisplaybreaks[4]

\theoremstyle{plain}
\newtheorem{thm}{Theorem}[section]
\newtheorem{lemm}[thm]{Lemma}

\theoremstyle{definition}
\newtheorem{df}[thm]{Definition}
\newtheorem{rem}[thm]{Remark}

\makeatletter

\@addtoreset{equation}{section}
\makeatother

\renewcommand{\div}{\operatorname{div}}
\newcommand{\dB}{\dot{B}}

\newcommand{\supp}{\operatorname{supp}}

\renewcommand{\leq}{\leqslant}
\renewcommand{\geq}{\geqslant}

\newcommand{\pnabla}{{\nabla}^{\perp}}

\newcommand{\n}[1]{{\left\|#1\right\|}}

\newcommand{\lp}[1]{\left[#1\right]}
\newcommand{\Mp}[1]{\left\{#1\right\}}
\renewcommand{\sp}[1]{\left(#1\right)}

\begin{document}
\title[stationary quasi-geostrophic equation]
{Sharp well-posedness and ill-posedness \\ of the stationary quasi-geostrophic equation}
\author[M.~Fujii]{Mikihiro Fujii}
\address[M.~Fujii]{Graduate School of Science, Nagoya City University, Nagoya, 467-8501, Japan}
\email[M.~Fujii]{fujii.mikihiro@nsc.nagoya-cu.ac.jp}
\author[T.~Iwabuchi]{Tsukasa Iwabuchi}
\address[T.~Iwabuchi]{
Mathematical Institute, Tohoku University, Sendai, 980-8578, Japan}
\email[T.~Iwabuchi]{t-iwabuchi@tohoku.ac.jp}
\begin{abstract}
We consider the stationary problem for the quasi-geostrophic equation on the whole plane and investigate its well-posedness and ill-posedness.
In \cite{Fuj-24}, it was shown that the two-dimensional stationary Navier--Stokes equations are ill-posed in the critical Besov spaces $\dB_{p,1}^{\frac{2}{p}-1}(\mathbb{R}^2)$ with $1 \leq p \leq 2$.
Although the quasi-geostrophic equation has the same invariant scale structure as the Navier--Stokes equations, we reveal that the quasi-geostrophic equation is well-posed in the scaling critical Besov spaces $\dB_{p,q}^{\frac{2}{p}-1}(\mathbb{R}^2)$ with $(p,q) \in [1,4) \times [1,\infty]$ or $(p,q)=(4,2)$ due to the better properties of the nonlinear structure of the quasi-geostrophic equation compared to that of the Navier--Stokes equations.
Moreover, we also prove the optimality for the above range of $(p,q)$ ensuring the well-posedness in the sense that the stationary quasi-geostrophic equation is ill-posed for all the other cases.
\end{abstract}
\maketitle

\section{Introduction}\label{sec:intro}
We consider the stationary problem of the quasi-geostrophic equation on the whole plane $\mathbb{R}^2$:
\begin{align}\label{eq:QG}
    \begin{cases}
        -\Delta \theta + u \cdot \nabla \theta = f, & \qquad x \in \mathbb{R}^2,\\ 
        u=\pnabla (-\Delta)^{-\frac{1}{2}} \theta, & \qquad x \in \mathbb{R}^2,
    \end{cases}
\end{align}
where $\pnabla=(-\partial_{x_2},\partial_{x_1})$.
Here, $\theta=\theta(x)$ and $u=(u_1(x),u_2(x))$ are the unknown potential temperature and velocity field of the fluid, respectively,
while $f=f(x)$ is the given external force.
The stationary quasi-geostrophic equation has the scaling invariant structure.
If $\theta$ is a solution to \eqref{eq:QG} for some given external force $f$, then $\theta_{\lambda}(x)=\lambda \theta(\lambda x)$ also solves \eqref{eq:QG} with $f_{\lambda}(x)=\lambda^3 f(\lambda x)$ for all $\lambda > 0$. We say that the data space $D$ and the solution space $S$ is scaling critical if $\n{f_{\lambda}}_D=\n{f}_D$ and $\n{\theta_{\lambda}}_S=\n{\theta}_S$ for all $\lambda>0$.
For instance, the homogeneous Besov spaces $D=\dB_{p,q}^{\frac{2}{p}-3}(\mathbb{R}^2)$ and $S=\dB_{p,q}^{\frac{2}{p}-1}(\mathbb{R}^2)$ ($1 \leq p,q \leq \infty$) are scaling critical.
The aim of this paper is to investigate the well-posedness and ill-posedness of the equation \eqref{eq:QG} in the scaling critical Besov spaces and provide the sharp result in the sense that it completely classifies whether \eqref{eq:QG} is well-posed or ill-posed for all $1 \leq p,q \leq \infty$.


\subsection*{Known results and the position of our study} 
Before stating the main results precisely, 
we recall the previous studies related to our work and mention the position of our study.
The quasi-geostrophic equation has received much attention for the case of the non-stationary initial value problem
\begin{equation}\label{eq:n-QG}
	\begin{cases}
		\partial_t\theta+(-\Delta)^{\alpha}\theta+u\cdot \nabla \theta=0,
			\quad & t>0,x\in \mathbb{R}^2,\\
		u=\pnabla (-\Delta)^{-\frac{1}{2}} \theta,
		 	\quad & t\geqslant 0,x\in \mathbb{R}^2,\\
		\theta(0,x)=\theta_0(x),
			\quad & x\in \mathbb{R}^2,
	\end{cases}
\end{equation}
where $0<\alpha \leq1$ and $\theta_0=\theta_0(x)$ is the given initial datum.
The local and global well-posedness for \eqref{eq:n-QG} are well-known by many researchers. 
See 
\cites{CC-04,CW-99,Iwa-Ued,Wu2,Wu} for the sub-critical case $\alpha>1/2$, 
\cites{Cha-Lee-03, CCW-01,DD-08, Kis-Naz-Vol-07,Zha-07} for the critical case $\alpha=1/2$, 
and 
\cites{Bae,Cha-Lee-03, Miu-06, Che-Mia-Zha-07} for the super-critical case $0< \alpha <1/2$.
Here, we remark that for the sub-critical and critical case $\alpha \geq 1/2$, \cites{CW-99,Kis-Naz-Vol-07} proved the global well-posedness for arbitrary large initial data in the scaling critical framework, whereas the global regularity for the super-critical case remains open.
Next, we mention the known results for the stationary quasi-geostrophic equation with the fractional Laplacian:
\begin{align}\label{sQG-a}
    \begin{cases}
        (-\Delta)^{\alpha} \theta + u \cdot \nabla \theta = f, \quad & x \in \mathbb{R}^2, \\
        u=\pnabla (-\Delta)^{-\frac{1}{2}} \theta,
		\quad & x\in \mathbb{R}^2.
    \end{cases}
\end{align}
There seem less results on \eqref{sQG-a} in comparison with the initial value problem on the non-stationary quasi-geostrophic equation.
For \eqref{sQG-a} with $\alpha=1/2$, Dai \cite{Dai-15} showed the existence of weak solutions $\theta \in H^{\frac{1}{2}}(\mathbb{R}^2)$ for small data $f \in W^{\frac{1}{2},4}(\mathbb{R}^2) \cap L^{\infty}(\mathbb{R}^2)$ satisfying $\widehat{f}(\xi)=0$ in some low frequency region. 
Moreover, she removed the condition $\widehat{f}(\xi)=0$ around $\xi = 0$ in the case of $1/2 < \alpha < 1$ and proved the existence of weak solutions $\theta \in H^{\alpha}(\mathbb{R}^2)$ for small data $f \in W^{1-\alpha,4}(\mathbb{R}^2) \cap L^{\frac{4}{2\alpha-1}}(\mathbb{R}^2)$.
For the unique existence of strong solutions in the sub-critical case $1/2 < \alpha < 1$,
Hadadifardo--Stefanov \cite{HS-21} proved that for any small external force $f \in L^{\frac{2}{2\alpha-1}}(\mathbb{R}^2)$, there exists a unique small solution $\theta \in \dot{W}^{-2\alpha,\frac{2}{2\alpha-1}}(\mathbb{R}^2)$; note that their framework is scaling critical.
Recently, Kozono--Kunstmann--Shimizu \cite{KKS} considered the generalized quasi-geostrophic equation on $\mathbb{R}^n$ with $n\geq 2$:
\begin{align}\label{eq:gQG}
    \begin{cases}
        (-\Delta)^{\alpha}\theta + u \cdot \nabla \theta = f, & \qquad x \in \mathbb{R}^n, \\
        u=\mathbb{P}S(\nabla(-\Delta)^{-\frac{1}{2}}\theta), & \qquad x \in \mathbb{R}^n,
    \end{cases}
\end{align}
where $\alpha>0$ and $S \in \mathbb{R}^{n \times n}$ are constant, and $\mathbb{P}=I+\nabla(-\Delta)^{-1}\div$ denotes the Helmholtz projection to the divergence free vector fields.
In \cite{KKS}, it was shown that if $1/2<\alpha <1/2 + n/4$, then for any $1 \leq p < n/(2\alpha -1)$ and $1 \leq q \leq \infty$, the small solution $\theta$ in the critical Besov space $\dB_{p,q}^{\frac{n}{p}+1-2\alpha}(\mathbb{R}^n)$ is constructed for small external force $f \in \dB_{p,q}^{\frac{n}{p}+1-4\alpha}(\mathbb{R}^n)$. 
Note that \eqref{eq:QG} corresponds to the case $n=2$ and $\alpha=1$, where \cite{KKS} does not encompass, and it was conjectured in \cite{KKS}*{Remarks (iii) after Theorem 2.3} that \eqref{eq:QG} is not solvable.
Indeed, the equation \eqref{eq:QG} has the same scaling structure as for the two-dimensional stationary Navier--Stokes equations, which is the following system with $n=2$:
\begin{align}\label{eq:sNS}
    \begin{cases}
        -\Delta u + \mathbb{P}(u \cdot \nabla) u = \mathbb{P} f, &\qquad x \in \mathbb{R}^n,\\
        \div u = 0, &\qquad x \in \mathbb{R}^n
    \end{cases}
\end{align} 
and the first author \cite{Fuj-24} proved the ill-posedness of \eqref{eq:sNS} with $n=2$ from $\dB_{p,1}^{\frac{2}{p}-3}(\mathbb{R}^2)$ to $\dB_{p,1}^{\frac{2}{p}-1}(\mathbb{R}^2)$ with $1 \leq p \leq 2$, which is the framework that \eqref{eq:sNS} is solvable in the higher dimensional case $n\geq 3$; see \cites{Che-93,Koz-Yam-95-PJA,Kan-Koz-Shi-19,Tsu-19-JMAA,Tsu-19-ARMA,Li-Yu-Zhu,Tan-Tsu-Zha-pre} for the details.
We notice that the unsolvability of \eqref{eq:sNS} with $n=2$ comes from fact that the key bilinear estimate 
\begin{align}
    \n{(-\Delta)^{-1}\mathbb{P}(u \cdot \nabla)v}_{\dB_{p,q}^{\frac{2}{p}-1}}
    \leq
    C
    \n{u}_{\dB_{p,q}^{\frac{2}{p}-1}}
    \n{v}_{\dB_{p,q}^{\frac{2}{p}-1}}
\end{align}
fails in almost all $1 \leq p,q \leq \infty$; see \cite{Fuj-pre} for the detail.
Despite such situation, in the present paper, we reveal that the nonlinear term $u \cdot \nabla \theta$ of the quasi-geostrophic equation \eqref{eq:QG} has a better regularity property than the two-dimensional Navier--Stokes equations and show that \eqref{eq:QG} is well-posed from $\dB_{p,q}^{\frac{2}{p}-3}(\mathbb{R}^2)$ to $\dB_{p,q}^{\frac{2}{p}-1}(\mathbb{R}^2)$ with $(p,q) \in ([1,4) \times [1,\infty]) \cup \{(4,2)\}$. 
Moreover, we show that the range $(p,q)\in ([1,4) \times [1,\infty]) \cup \{(4,2)\}$ that ensures the well-posed is sharp by showing the ill-posedness for the other case of $(p,q)$.


\subsection*{Reformulation of the problem and the statements of our main results}
Following the argument in \cites{Kan-Koz-Shi-19,KKS},
we rewrite \eqref{eq:QG} as  
\begin{align}\label{eq:rQG}
    \theta = \mathcal{L}f + (-\Delta)^{-1}\div (\theta u),\qquad
    u=\pnabla (-\Delta)^{-\frac{1}{2}}\theta,
\end{align}
where 
\begin{align}\label{L}
    \mathcal{L}f:=(-\Delta)^{-1}f,
\end{align}
and
try to control the product term $(-\Delta)^{-1}\div (\theta_1 u_2)$ with $u_2=\nabla^{\perp}(-\Delta)^{-\frac{1}{2}}\theta_2$ as 
\begin{align}
    \n{(-\Delta)^{-1}\div(\theta_1 u_2)}_{\dB_{p,q}^{\frac{2}{p}-1}}
    &=
    \n{(-\Delta)^{-1}\div\sp{\theta_1 \nabla^{\perp}(-\Delta)^{-\frac{1}{2}}\theta_2}}_{\dB_{p,q}^{\frac{2}{p}-1}}\\
    &\leq 
    C
    \n{\theta_1 \nabla^{\perp}(-\Delta)^{-\frac{1}{2}}\theta_2}_{\dB_{p,q}^{\frac{2}{p}-2}}.
\end{align}
However, we may not proceed the above estimate since the para-product estimate
\begin{align}
    \n{fg}_{\dB_{p,q}^{\frac{2}{p}-2}}
    \leq 
    C
    \n{f}_{\dB_{p,q}^{\frac{2}{p}-1}}
    \n{g}_{\dB_{p,q}^{\frac{2}{p}-1}}
\end{align}
{\it fails} for all $1 \leq p,q \leq \infty$.
Due to the similar circumstance, the two-dimensional Navier--Stokes equations \eqref{eq:sNS} with $n=2$ is ill-posed even in the narrowest critical Besov space and \cite{KKS} did not treat the case $n=2$ and $\alpha=1$ in \eqref{eq:gQG}.

In the present paper, we reveal that, in contrast to the Navier--Stokes equations \eqref{eq:sNS}, the nonlinear term of the quasi-geostrophic equation \eqref{eq:QG} possesses a better nonlinear structure.
To see this, we rewrite the nonlinear term of \eqref{eq:rQG} as 
\begin{align}
    \mathscr{F}\lp{(-\Delta)^{-1}\div (\theta u)}(\xi)
    ={}&
    \mathscr{F}\lp{(-\Delta)^{-1}\div \sp{\theta \pnabla(-\Delta)^{-\frac{1}{2}}\theta}}(\xi)
    \\
    ={}&
    -
    \int_{\mathbb{R}^2}
    \frac{\xi\cdot \eta^{\perp}}{|\xi|^2|\eta|}
    \widehat{\theta}(\xi-\eta)
    \widehat{\theta}(\eta)
    d\eta.
\end{align}
Here, using 
\begin{align}
    -
    \int_{\mathbb{R}^2}
    \frac{\xi\cdot \eta^{\perp}}{|\xi|^2|\eta|}
    \widehat{\theta}(\xi-\eta)
    \widehat{\theta}(\eta)
    d\eta
    &=
    -
    \int_{\mathbb{R}^2}
    \frac{\xi \cdot (\xi-\eta)^{\perp}}{|\xi|^2|\xi-\eta|}
    \widehat{\theta}(\eta)
    \widehat{\theta}(\xi-\eta)
    d\eta\\
    &=
    \int_{\mathbb{R}^2}
    \frac{\xi\cdot \eta^{\perp}}{|\xi|^2|\xi-\eta|}
    \widehat{\theta}(\xi-\eta)
    \widehat{\theta}(\eta)
    d\eta,
\end{align}
we have
\begin{align}
    \mathscr{F}\lp{(-\Delta)^{-1}\div (\theta u)}(\xi)
    ={}&
    \frac{1}{2}
    \int_{\mathbb{R}^2}
    \frac{\xi\cdot \eta^{\perp}}{|\xi|^2}
    \sp{\frac{1}{|\xi-\eta|}-\frac{1}{|\eta|}}
    \widehat{\theta}(\xi-\eta)
    \widehat{\theta}(\eta)
    d\eta    \\
    ={}&
    \frac{1}{2}
    \int_{\mathbb{R}^2}
    \frac{\xi\cdot \eta^{\perp}}{|\xi|^2}
    \frac{\xi\cdot(2\eta-\xi)}{|\eta|+|\xi-\eta|}
    \frac{\widehat{\theta}(\xi-\eta)}{|\xi-\eta|}
    \frac{\widehat{\theta}(\eta)}{|\eta|}
    d\eta.
\end{align}
Therefore, setting 
\begin{align}
    \mathcal{B}[\theta_1, \theta_2]
    :={}&
    (-\Delta)^{-1} \div \div 
    \mathcal{T}
    \lp{(-\Delta)^{-\frac{1}{2}}\theta_1, \pnabla(-\Delta)^{-\frac{1}{2}}\theta_2}, 
    \label{B}
    \\
    \mathcal{T}[\theta, v]
    :={}&
    \mathscr{F}^{-1}
    \lp{
    \int_{\mathbb{R}^2}
    \frac{\xi-2\eta}{2(|\eta|+|\xi-\eta|)}
    \otimes 
    \sp{\widehat{\theta}(\xi-\eta)\widehat{v}(\eta)}
    d\eta
    }\label{T}
\end{align}
for appropriate given functions $\theta_1,\theta_2,\theta:\mathbb{R}^2 \to \mathbb{R}$ and $v:\mathbb{R}^2 \to \mathbb{R}^2$,
we have 
\begin{align}
    (-\Delta)^{-1}\div (u\theta) = \mathcal{B}[\theta,\theta],\qquad u = \pnabla(-\Delta)^{-\frac{1}{2}}\theta.
\end{align}
In comparison to the single divergence form $(-\Delta)^{-1}\div (u\theta)$, the rewritten formula $\mathcal{B}[\theta,\theta]$ has the double divergence form, which help us from the difficulty stated above and succeed in closing the nonlinear estimates in some scaling critical Besov spaces; see Lemma \ref{lemm:nonlin-1} below.
Motivated by this observation, we define the notion of the solution to \eqref{eq:QG} as follows.
\begin{df}
    For $1 \leq p,q \leq \infty$ and $f \in \dB_{p,q}^{\frac{2}{p}-3}(\mathbb{R}^2)$,
    it is said that $\theta \in \dB_{p,q}^{\frac{2}{p}-1}(\mathbb{R}^2)$ is a solution to \eqref{eq:QG} if it satisfies
    \begin{align}
        \theta 
        =
        \mathcal{L}f
        +
        \mathcal{B}[\theta,\theta],
    \end{align}
    where the operators $\mathcal{L}$ and $\mathcal{B}[\cdot,\cdot]$ are defined in \eqref{L} and \eqref{B}, respectively.
\end{df}

The following theorem is the first main result of this paper that states the well-posedness of \eqref{eq:QG} in the scaling critical Besov spaces.
\begin{thm}[Well-posedness of \eqref{eq:QG}]\label{thm:WP}
    Let $1 \leq p < 4$ and $1 \leq q \leq \infty$, or $p=4$ and $q =2$.
    Then, there exist positive constants $\delta_0=\delta_0(p,q)$ and $\varepsilon_0=\varepsilon_0(p,q)$ such that 
    if the external force $f \in D_{p,q}(\mathbb{R}^2)$, where
    \begin{align}
        D_{p,q}(\mathbb{R}^2):=\Mp{f \in \dB_{p,q}^{\frac{2}{p}-3}(\mathbb{R}^2)\ ;\ \n{f}_{\dB_{p,q}^{\frac{2}{p}-3}} \leq \delta_0},
    \end{align}
    then there exists a unique solution $\theta \in S_{p,q}(\mathbb{R}^2)$, where 
    \begin{align}
        S_{p,q}(\mathbb{R}^2):=\Mp{
        \theta \in \dB_{p,q}^{\frac{2}{p}-1}(\mathbb{R}^2)\ ;\ 
        \n{\theta}_{\dB_{p,q}^{\frac{2}{p}-1}}\leq \varepsilon_0
        }.
    \end{align}
    Moreover, the solution map 
    $\Phi:D_{p,q}(\mathbb{R}^2) \ni f \mapsto \theta \in S_{p,q}(\mathbb{R}^2)$
    is Lipschitz continuous
    from 
    $\dB_{p,q}^{\frac{2}{p}-3}(\mathbb{R}^2)$
    to
    $\dB_{p,q}^{\frac{2}{p}-1}(\mathbb{R}^2)$.
\end{thm}

Next, we focus on the optimality of the range of $(p,q)$ for well-posedness and state that \eqref{eq:QG} is ill-posed for indexes other than those whose well-posedness is guaranteed by the aforementioned theorem. 
\begin{thm}[Ill-posedness of \eqref{eq:QG}]\label{thm:IP}
    Let $p > 4$ and $1 \leq q \leq \infty$, or $p=4$ and $1 \leq q \leq \infty$ with $q \neq 2$.
    Then, \eqref{eq:QG} is ill-posed
    in the sense that 
    the solution map 
    $f \mapsto \theta$
    is discontinuous
    from 
    $\dB_{p,q}^{\frac{2}{p}-3}(\mathbb{R}^2)$
    to
    $\dB_{p,q}^{\frac{2}{p}-1}(\mathbb{R}^2)$ at the origin.
    More precisely, the following statements hold true.
    \begin{enumerate}
        \item 
        For $4 < p\leq \infty$ and $1 \leq q \leq \infty$, or $p=4$ and $2 <q \leq \infty$, there exists a sequence $\{f_N\}_{N=1}^{\infty} \subset \mathscr{S}(\mathbb{R}^2)$ of the external forces satisfying 
        \begin{align}
            \lim_{N \to \infty}
            \n{f_N}_{\dB_{p,q}^{\frac{2}{p}-3}} = 0,
        \end{align}
        whereas the corresponding sequence $\{\theta_N\}_{N=1}^{\infty} \subset \dB_{p,q}^{\frac{2}{p}-1}(\mathbb{R}^2)$ of solutions fulfills
        \begin{align}
            \liminf_{N\to \infty}
            \n{\theta_N}_{\dB_{\infty,\infty}^{-1}}
            >
            0.
        \end{align}
        \item 
        For $p=4$ and $1 \leq q < 2$, there exists a sequence $\{f_N\}_{N=1}^{\infty} \subset \mathscr{S}(\mathbb{R}^2)$ of the external forces such that
        \begin{align}
            \lim_{N \to \infty}
            \n{f_N}_{\dB_{4,q}^{-\frac{5}{2}}} = 0,
        \end{align}
        while the corresponding sequence $\{\theta_N\}_{N=1}^{\infty} \subset \dB_{4,q}^{-\frac{1}{2}}(\mathbb{R}^2)$ of solutions satisfies
        \begin{align}
            \lim_{N\to \infty}
            \n{\theta_N}_{\dB_{4,q}^{-\frac{1}{2}}}
            =
            \infty
        \end{align}
    \end{enumerate}
\end{thm}
\begin{rem}
We give some comments on Theorem \ref{thm:IP}.
For the case of $(p,q) \in ((4,\infty) \times [1,\infty])\cup(\{4\} \times (2,\infty))$, we actually show that the continuity of the solution map fails in weaker topology $\dB_{p,q}^{\frac{2}{p}-3}(\mathbb{R}^2) \to \dB_{\infty, \infty}^{-1}(\mathbb{R}^2)$, which is similar to \cite{Tsu-19-ARMA}.
For the case of $(p,q) \in \{ 4 \} \times [1,2)$, we succeed in proving the norm-inflation phenomenon $\n{\theta_N }_{\dB_{4,q}^{-\frac{1}{2}}}\geq cN^{\frac{1}{q}-\frac{1}{2}}/(\log N)^2$ that has not been proved for the ill-posedness of the stationary Navier--Stokes equations; see \cites{Fuj-24,Tsu-19-ARMA,Li-Yu-Zhu,Tan-Tsu-Zha-pre}.
\end{rem}
This paper is organized as follows.
In Section \ref{sec:pre}, we recall the notations and definition of the Besov spaces.
In Sections \ref{sec:pf_p<4} and \ref{sec:pf_p>4}, we provide the proof of Theorems \ref{thm:WP} and \ref{thm:IP}, respectively.
\section{Preliminaries}\label{sec:pre}
In this section, we introduce notations and function spaces, which are to be used in this paper.
Throughout this paper, we denote by $C\geq 1$ and $0<c<1$ the constants, which may differ in each line. 
In particular, $C=C(a_1,...,a_n)$ means that $C$ depends only on $a_1,...,a_n$.
For two non-negative numbers $A$, $B$, the relation $A \sim B$ means that there exists positive constant $C$ such that $C^{-1}A\leqslant B \leqslant C A$ holds.
Let $\mathscr{S}(\mathbb{R}^2)$ be the set of all Schwartz functions on $\mathbb{R}^2$ and $\mathscr{S}'(\mathbb{R}^2)$ represents the 
set of all tempered distributions on $\mathbb{R}^2$.
We use $L^p(\mathbb{R}^2)$  with $1 \leqslant p \leqslant \infty$ to denote the standard Lebesgue spaces on $\mathbb{R}^2$.
For $f \in \mathscr{S}(\mathbb{R}^2)$, the Fourier transform and inverse Fourier transform of $f$ are defined as
\begin{align}
    \mathscr{F}[f](\xi)=\widehat{f}(\xi):=\int_{\mathbb{R}^2} e^{-ix\cdot\xi} f(x)dx,\qquad
    \mathscr{F}^{-1}[f](x):=\frac{1}{(2\pi)^2}\int_{\mathbb{R}^2} e^{ix\cdot\xi} f(\xi)d\xi.
\end{align}

Let $\{\phi_j\}_{j \in \mathbb{Z}} \subset \mathscr{S}(\mathbb{R}^2)$ be a dyadic partition of unity satisfying 
\begin{align}
    &
    0 \leqslant \widehat{\phi_0}(\xi) \leqslant 1,\\
    &
    \supp \widehat{\phi_0} \subset \{ \xi \in \mathbb{R}^2\ ;\ 2^{-1} \leqslant |\xi| \leqslant 2 \},\\
    &
    \widehat{\phi_0}(\xi) = 1 \qquad {\rm for\ all\ }\xi \in \mathbb{R}^2{\rm\ with\ }\frac{7}{8} \leqslant |\xi| \leqslant \frac{5}{4},
\end{align}
and
\begin{align}
    \sum_{j \in \mathbb{Z}}
    \widehat{\phi_j}(\xi)
    =1
    \qquad {\rm for\ all\ }\xi \in \mathbb{R}^2 \setminus \{ 0 \},
\end{align}
where $\widehat{\phi_j}(\xi) = \widehat{\phi_0}(2^{-j}\xi)$.
We then define the homogeneous Besov spaces $\dB_{p,q}^s(\mathbb{R}^2)$ ($1 \leqslant p,q \leqslant \infty$, $s \in \mathbb{R}$) by 
\begin{align}
    \dB_{p,q}^s(\mathbb{R}^2)
    :={}&
    \left\{
    f \in \mathscr{S}'(\mathbb{R}^2) / \mathscr{P}(\mathbb{R}^2)
    \ ; \ 
    \| f \|_{\dB_{p,q}^s}
    <
    \infty
    \right\},\\
    \| f \|_{\dB_{p,q}^s}
    :={}&
    \left\|
    \left\{
    2^{sj}
    \| \phi_j * f \|_{L^p}
    \right\}_{j \in \mathbb{Z}}
    \right\|_{\ell^{q}},
\end{align}
where $\mathscr{P}(\mathbb{R}^2)$ denotes the set of all polynomials on $\mathbb{R}^2$.
It is well-known that if $s <2/p$ or $(s,q)=(2/p,1)$, then $\dB_{p,q}^s(\mathbb{R}^2)$ is identified as 
\begin{align}\label{chara-Besov}
    \dB_{p,q}^s(\mathbb{R}^2)
    \sim 
    \left\{
    f \in \mathscr{S}'(\mathbb{R}^2)\ ;\ 
    f = \sum_{j \in \mathbb{Z}} \phi_j * f
    \quad {\rm in\ }\mathscr{S}'(\mathbb{R}^2),\quad
    \| f \|_{\dB_{p,q}^s}
    <
    \infty
    \right\}.
\end{align}
See \cite{Saw-18}*{Theorem 2.31} for the proof of \eqref{chara-Besov}.
We refer to \cite{Saw-18} for the basic properties of Besov spaces.
\section{Proof of Theorem \ref{thm:WP}}\label{sec:pf_p<4}
In this section, we present the proof of Theorem \ref{thm:WP}.
\subsection*{Key lemmas}
Before starting the proof of Theorem \ref{thm:WP}, 
we provide lemmas playing key roles in our calculus.
The following lemma plays the most crucial role for controlling the nonlinearity of \eqref{eq:QG}.
\begin{lemm}\label{lemm:nonlin-1}
    Let 
    $1 \leq p < 4$ and $1 \leq q \leq \infty$, 
    or 
    $p=4$ and $q=2$.
    Then, there exists a positive constant $C=C(p,q)$ such that 
    \begin{align}\label{nonlin-1}
        \n{\mathcal{B}[f,g]}_{\dB_{p,q}^{\frac{2}{p}-1}}
        \leq
        C
        \n{f}_{\dB_{p,q}^{\frac{2}{p}-1}}
        \n{g}_{\dB_{p,q}^{\frac{2}{p}-1}}
    \end{align}
    for all $f,g\in \dB_{p,q}^{\frac{2}{p}-1}(\mathbb{R}^2)$.
\end{lemm}
To prove Lemma \ref{lemm:nonlin-1}, we prepare a fact for the boundedness of bilinear Fourier multiplier:
\begin{lemm}[\cite{Gra-Miy-Tom-13}]\label{lemm:GMT}
    Let $m \in C^2((\mathbb{R}^2 \setminus \{ 0\}) \times (\mathbb{R}^2 \setminus \{ 0\}))$ satisfy
    \begin{align}
        \left| \partial_{\xi}^{\alpha} \partial_{\eta}^{\beta} m(\xi,\eta) \right|
        \leq
        C_{\alpha,\beta}
        \sp{|\xi|+|\eta|}^{-(|\alpha|+|\beta|)}, \qquad 
        \xi,\eta \in \mathbb{R}^2
    \end{align}
    for $\alpha, \beta \in (\mathbb{N}\cup\{0\})^2$ with $|\alpha|,|\beta| \leq 2$ and for some positive constant $C_{\alpha,\beta}$.
    Let $T_m$ be a bilinear operator defined by 
    \begin{align}
        T_m[f,g](x)
        :={}&
        \mathscr{F}^{-1}
        \lp{\int_{\mathbb{R}^2}m(\xi-\eta,\eta)\widehat{f}(\xi-\eta)\widehat{g}(\eta)d\eta}(x)\\
        ={}&
        \frac{1}{(2\pi)^2}
        \int_{\mathbb{R}^2_{\xi} \times \mathbb{R}^2_{\eta}}
        e^{ix\cdot(\xi+\eta)}
        m(\xi,\eta)\widehat{f}(\xi)\widehat{g}(\eta)d\xi d\eta.
    \end{align}
    Let $p,p_1,p_2$ satisfy $1/p = 1/p_1 + 1/p_2$. 
    Assume $1 < p \leq 2$ and $1 < p_1,p_2 \leq\infty$, or $2 \leq p,p_1,p_2 < \infty$.
    Then,
    there exists a positive constant $C=C(p,p_1,p_2)$ such that 
    \begin{align}
        \n{T_m[f,g]}_{L^p}
        \leq
        C\n{f}_{L^{p_1}}\n{g}_{L^{p_2}}
    \end{align}
    for all $f \in L^{p_1}(\mathbb{R}^2)$ and $g \in L^{p_2}(\mathbb{R}^2)$.
    Moreover, there holds
    \begin{align}
        \n{T_m[f,g]}_{L^1}
        \leq
        C\n{f}_{L^{\infty}}\n{g}_{\mathcal{H}^1}
    \end{align}
    for all $f \in L^{\infty}(\mathbb{R}^2)$ and $g \in \mathcal{H}^1(\mathbb{R}^2)$.
    Here, $\mathcal{H}^1(\mathbb{R}^2)$ denotes the Hardy space on $\mathbb{R}^2$ with the integrability index $1$.
\end{lemm}
We are in a position to show Lemma \ref{lemm:nonlin-1}.
\begin{proof}
    Recalling the Bony decomposition, 
    we decompose $\mathcal{B}[f,g]$ as
    \begin{align}
        \mathcal{B}[f,g]
        =\mathcal{B}_1[f,g]+\mathcal{B}_2[f,g]+\mathcal{B}_3[f,g],
    \end{align}
    where
    \begin{align}
        \mathcal{B}_1[f,g]
        :={}&
        \sum_{k \leq \ell -3}
        \mathcal{B}\lp{f_k,g_{\ell}}\\
        ={}&
        \frac{1}{2}
        \sum_{k \leq \ell -3}
        (-\Delta)^{-1}\div
        \Mp{f_k\sp{\pnabla(-\Delta)^{-\frac{1}{2}}g_{\ell}} + \sp{\pnabla(-\Delta)^{-\frac{1}{2}}f_k}g_{\ell}},\\
        \mathcal{B}_2[f,g]
        :={}&
        \mathcal{B}_1[g,f],\\
        \mathcal{B}_3[f,g]
        :={}&
        \sum_{|k - \ell|\leq 2}
        \mathcal{B}\lp{f_k,g_{\ell}}.
    \end{align}
    Here, we have used the abbreviation $h_k:=\phi_k * h$ for $k \in \mathbb{Z}$ and $h \in \mathscr{S}'(\mathbb{R}^2)/\mathscr{P}(\mathbb{R}^2)$, and also used the following identity
    \begin{align}
        \mathcal{B}[f_k,g_\ell]
        =
        \frac{1}{2}
        (-\Delta)^{-1}\div
        \Mp{f_k\sp{\pnabla(-\Delta)^{-\frac{1}{2}}g_{\ell}} + \sp{\pnabla(-\Delta)^{-\frac{1}{2}}f_k}g_{\ell}},
    \end{align}
    which is obtained by the similar calculations as in the motivation of the definition of $\mathcal{B}[\cdot,\cdot]$ in Section \ref{sec:intro}.
    By the standard paraproduct estimates, there holds
    \begin{align}
        &\n{\mathcal{B}_1[f,g]}_{\dB_{p,q}^{\frac{2}{p}-1}}\\
        &\quad
        \leq{}
        C
        \n{
        \Mp{
        2^{(\frac{2}{p}-2)j}
        \sum_{|m|\leq 3}
        \sum_{k \leq j+m -3}
        \n{f_k}_{L^{\infty}}
        \n{\pnabla(-\Delta)^{-\frac{1}{2}}g_{j+m}}_{L^p}
        }_{j \in \mathbb{Z}}
        }_{\ell^{q}}\\
        &
        \qquad
        +
        C
        \n{
        \Mp{
        2^{(\frac{2}{p}-2)j}
        \sum_{|m|\leq 3}
        \sum_{k \leq j+m -3}
        \n{\pnabla(-\Delta)^{-\frac{1}{2}}f_k}_{L^{\infty}}\n{g_{j+m}}_{L^p}
        }_{j \in \mathbb{Z}}
        }_{\ell^{q}}\\
        &\quad
        \leq{}
        C
        \n{
        \Mp{
        \sum_{|m|\leq 3}
        \Mp{\sum_{k \leq j+m -3}
        \sp{2^{-k}\n{f_k}_{L^{\infty}}}^q}^{\frac{1}{q}}
        2^{(\frac{2}{p}-1)(j+m)}\n{g_{j+m}}_{L^p}
        }_{j \in \mathbb{Z}}
        }_{\ell^{q}}\\
        &\quad 
        \leq{}
        C
        \n{f}_{\dB_{\infty,q}^{-1}}
        \n{g}_{\dB_{p,q}^{\frac{2}{p}-1}}
        \leq
        C
        \n{f}_{\dB_{p,q}^{\frac{2}{p}-1}}
        \n{g}_{\dB_{p,q}^{\frac{2}{p}-1}}.
    \end{align}
    Similarly, we have
    \begin{align}
        \n{\mathcal{B}_2[f,g]}_{\dB_{p,q}^{\frac{2}{p}-2}}
        \leq
        C
        \n{f}_{\dB_{p,q}^{\frac{2}{p}-1}}
        \n{g}_{\dB_{p,q}^{\frac{2}{p}-1}}.
    \end{align}
    We note that the estimates for $\mathcal{B}_1[f,g]$ and $\mathcal{B}_2[f,g]$ hold for all $1 \leq p,q \leq \infty$.
    For the estimate of $\mathcal{B}_3[f,g]$, we see that 
    \begin{align}
        \n{\mathcal{B}_3[f,g]}_{\dB_{p,q}^{\frac{2}{p}-2}}
        \leq{}&
        C
        \n{\sum_{|k-\ell|\leq 2}\mathcal{T}\lp{{(-\Delta)^{-\frac{1}{2}}f_k},{\nabla^{\perp}(-\Delta)^{-\frac{1}{2}}g_{\ell}}}}_{\dB_{p,q}^{\frac{2}{p}-1}},
    \end{align}
    where $\mathcal{T}[\cdot,\cdot]$ is defined in \eqref{T}.
    We will apply Lemma \ref{lemm:GMT} to control the bilinear operator $\mathcal{T}[\cdot,\cdot]$.
    For the case of $(p,q)=(4,2)$, it holds
    \begin{align}
        \n{\mathcal{B}_3[f,g]}_{\dB_{4,2}^{-\frac{1}{2}}}
        \leq{}&
        C
        \n{
        \sum_{|k-\ell|\leq 2}
        \mathcal{T}\lp{{(-\Delta)^{-\frac{1}{2}}f_k},{\nabla^{\perp}(-\Delta)^{-\frac{1}{2}}g_{\ell}}}
        }_{L^2}\\
        \leq{}&
        C
        \sum_{|k-\ell|\leq 2}
        \n{
        \mathcal{T}\lp{{(-\Delta)^{-\frac{1}{2}}f_k},{\nabla^{\perp}(-\Delta)^{-\frac{1}{2}}g_{\ell}}}}_{L^2}
        \\
        \leq{}&
        C
        \sum_{|k-\ell|\leq 2}
        \n{(-\Delta)^{-\frac{1}{2}}f_k}_{L^4}\n{\pnabla(-\Delta)^{-\frac{1}{2}}g_{\ell}}_{L^4}
        \\
        \leq{}&
        C
        \sum_{|k-\ell|\leq 2}
        2^{-k}\n{f_k}_{L^4}\n{g_{\ell}}_{L^4}
        \\
        \leq{}&
        C
        \n{f}_{\dB_{4,2}^{-\frac{1}{2}}}
        \n{g}_{\dB_{4,2}^{-\frac{1}{2}}}
    \end{align}
    For the case of $2 \leq p < 4$, we see that 
    \begin{align}
        \n{\mathcal{B}_3[f,g]}_{\dB_{p,q}^{\frac{2}{p}-1}}
        \leq{}&
        C
        \n{
        \sum_{|k-\ell|\leq 2}
        T\lp{{(-\Delta)^{-\frac{1}{2}}f_k},{\nabla^{\perp}(-\Delta)^{-\frac{1}{2}}g_{\ell}}}
        }_{\dB_{\frac{p}{2},q}^{\frac{4}{p}-1}}\\
        \leq{}&
        C
        \n{
        \Mp{
        2^{(\frac{4}{p}-1)j}
        \sum_{k \geq j-4}
        \sum_{|k-\ell|\leq 2}
        \n{
        T\lp{{(-\Delta)^{-\frac{1}{2}}f_k},{\nabla^{\perp}(-\Delta)^{-\frac{1}{2}}g_{\ell}}}}_{L^{\frac{p}{2}}}
        }_{j \in \mathbb{Z}}
        }_{\ell^q}\\
        \leq{}&
        C
        \n{
        \Mp{
        2^{(\frac{4}{p}-1)j}
        \sum_{k \geq j-4}
        \sum_{|k-\ell|\leq 2}
        \n{(-\Delta)^{-\frac{1}{2}}f_k}_{L^p}\n{g_{\ell}}_{L^p}
        }_{j \in \mathbb{Z}}
        }_{\ell^q}\\
        \leq{}&
        C
        \n{f}_{\dB_{p,q}^{\frac{2}{p}-1}}
        \n{g}_{\dB_{p,q}^{\frac{2}{p}-1}}
    \end{align}
    For the case of $1 < p < 2$, we have 
    \begin{align}
        \n{\mathcal{B}_3[f,g]}_{\dB_{p,q}^{\frac{2}{p}-1}}
        \leq{}&
        C
        \n{
        \Mp{
        2^{(\frac{2}{p}-1)j}
        \sum_{k \geq j-4}
        \sum_{|k-\ell|\leq 2}
        \n{
        T\lp{{(-\Delta)^{-\frac{1}{2}}f_k},{\nabla^{\perp}(-\Delta)^{-\frac{1}{2}}g_{\ell}}}}_{L^{p}}
        }_{j \in \mathbb{Z}}
        }_{\ell^q}\\
        \leq{}&
        C
        \n{
        \Mp{
        2^{(\frac{2}{p}-1)j}
        \sum_{k \geq j-4}
        \sum_{|k-\ell|\leq 2}
        \n{(-\Delta)^{-\frac{1}{2}}f_k}_{L^{\infty}}\n{g_{\ell}}_{L^p}
        }_{j \in \mathbb{Z}}
        }_{\ell^q}\\
        \leq{}&
        C
        \n{f}_{\dB_{\infty,q}^{-1}}
        \n{g}_{\dB_{p,q}^{\frac{2}{p}-1}}
        \leq{}
        C
        \n{f}_{\dB_{p,q}^{\frac{2}{p}-1}}
        \n{g}_{\dB_{p,q}^{\frac{2}{p}-1}}.
    \end{align}
    For the case of $p=1$, it follows that
    \begin{align}
        \n{\mathcal{B}_3[f,g]}_{\dB_{p,q}^{\frac{2}{p}-1}}
        \leq{}&
        C
        \n{
        \Mp{
        2^{j}
        \sum_{k \geq j-4}
        \sum_{|k-\ell|\leq 2}
        \n{
        T\lp{{(-\Delta)^{-\frac{1}{2}}f_k},{\nabla^{\perp}(-\Delta)^{-\frac{1}{2}}g_{\ell}}}}_{L^1}
        }_{j \in \mathbb{Z}}
        }_{\ell^q}\\
        \leq{}&
        C
        \n{
        \Mp{
        2^{j}
        \sum_{k \geq j-4}
        \sum_{|k-\ell|\leq 2}
        \n{(-\Delta)^{-\frac{1}{2}}f_k}_{L^{\infty}}\n{g_{\ell}}_{\mathcal{H}^1}
        }_{j \in \mathbb{Z}}
        }_{\ell^q}\\
        \leq{}&
        C
        \n{f}_{\dB_{\infty,q}^{-1}}
        \n{g}_{\dB_{1,q}^{1}}
        \leq{}
        C
        \n{f}_{\dB_{1,q}^{1}}
        \n{g}_{\dB_{1,q}^{1}},
    \end{align}
    where we have used  
    \begin{align}
        \n{
        \Mp{
        2^{\ell}
        \n{g_{\ell}}_{\mathcal{H}^1}
        }_{j \in \mathbb{Z}}
        }_{\ell^q}
        \sim
        \n{g}_{\dB_{1,q}^{1}},
    \end{align}
    which is implied by
    \begin{align}
        \n{g_{\ell}}_{\mathcal{H}^1}
        \sim 
        \n{g_{\ell}}_{\dot{F}_{1,2}^0}
        \sim 
        \n{g_{\ell-1}}_{L^1}
        +
        \n{g_{\ell}}_{L^1}
        +
        \n{g_{\ell+1}}_{L^1}.
    \end{align}
    Here, $\dot{F}_{1,2}^0(\mathbb{R}^2)$ denotes the Lizorkin--Tribel space.
    Thus, we complete the proof.
\end{proof}
\subsection*{Proof of Theorem \ref{thm:WP}}
Now, we are in a position to present the proof of Theorem \ref{thm:WP}.
\begin{proof}[Proof of Theorem \ref{thm:WP}]
Let $C_0=C_0(p)$ be a positive constant satisfying 
\begin{align}
    \n{\mathcal{L}f}_{\dB_{p,q}^{\frac{2}{p}-1}}
    \leq
    C_0
    \n{f}_{\dB_{p,q}^{\frac{2}{p}-3}}.
\end{align}
Let $C=C_1$ be the positive constant appearing in Lemma \ref{lemm:nonlin-1}.
Let $f \in \dB_{p,q}^{\frac{2}{p}-3}(\mathbb{R}^2)$ satisfy 
\begin{align}
    \n{f}_{\dB_{p,q}^{\frac{2}{p}-3}}
    \leq
    \frac{1}{8C_0C_1}.
\end{align}
Let us define the complete metric space $(X,d_X)$ by
\begin{align}
    &X:=
    \Mp{
    \theta \in \dB_{p,q}^{\frac{2}{p}-1}(\mathbb{R}^2)\ ; \ 
    \n{\theta}_{\dB_{p,q}^{\frac{2}{p}-1}}
    \leq
    \frac{1}{4C_1}.
    },\\
    &
    d_X(f,g):=\n{f-g}_{\dB_{p,q}^{\frac{2}{p}-1}}.
\end{align}
We consider a map $\Phi[\cdot]$ on $X$ given by
\begin{align}
    \Phi[\theta]:=\mathcal{L}f + \mathcal{B}[\theta,\theta], \qquad \theta \in X.
\end{align}
For each $\theta \in X$,
it follows from Lemma \ref{lemm:nonlin-1} that 
\begin{align}
    \n{\Phi[\theta]}_{\dB_{p,q}^{\frac{2}{p}-1}}
    \leq{}&
    C_0
    \n{f}_{\dB_{p,q}^{\frac{2}{p}-3}}
    +
    C_1
    \n{\theta}_{\dB_{p,q}^{\frac{2}{p}-1}}^2\\
    \leq{}&
    \frac{1}{8C_1}
    +
    \frac{1}{16C_1}
    \leq
    \frac{1}{4C_1},
\end{align}
which implies $\Phi[\theta] \in X$.
Let $\theta_1,\theta_2 \in X$.
Then, since
\begin{align}
    \Phi[\theta_1]-\Phi[\theta_2]
    =&
    \mathcal{B}[\theta_1-\theta_2,\theta_1]
    +
    \mathcal{B}[\theta_2,\theta_1-\theta_2],
\end{align}
we see that 
\begin{align}
    \n{\Phi[\theta_1]-\Phi[\theta_2]}_{\dB_{p,q}^{\frac{2}{p}-1}}
    \leq{}&
    C_1\sp{\n{\theta_1}_{\dB_{p,q}^{\frac{2}{p}-1}} + \n{\theta_2}_{\dB_{p,q}^{\frac{2}{p}-1}}}\n{\theta_1-\theta_2}_{\dB_{p,q}^{\frac{2}{p}-1}}\\
    \leq{}&
    \frac{1}{2}
    \n{\theta_1-\theta_2}_{\dB_{p,q}^{\frac{2}{p}-1}}.
\end{align}
Hence, $\Phi[\cdot]$ is a contraction map on $(X,d_X)$ and thus the Banach fixed point theorem yields the unique existence of $\theta \in X$ satisfying $\theta = \Phi[\theta]$.
This completes the proof.
\end{proof}

\section{Proof of Theorem \ref{thm:IP}}\label{sec:pf_p>4}
We are ready to prove Theorem \ref{thm:IP}.
\begin{proof}[Proof of Theorem \ref{thm:IP}]
Let $0 < \delta \leq 1/100$ be a positive constant to be determined depending only on $p$ and $q$ later.
In the following, we divide the proof into three steps:
\begin{itemize}
    \item [Step 1.] $4 < p \leq \infty$ and $1 \leq q \leq \infty$,
    \item [Step 2.] $p=4$ and $2 < q \leq \infty$,
    \item [Step 3.] $p=4$ and $1 \leq q < 2$.
\end{itemize}
We mention the outline of the argument in each step.
We first define the squeezing sequence $\{ f_N \}_{N=1}^{\infty} \subset \mathscr{S}(\mathbb{R}^2)$ of the external forces  and then establish the estimates for the first and second iterations
\begin{align}
    \theta_N^{(1)}
    :=\mathcal{L}f_N, \qquad
    \theta_N^{(2)}
    :=\mathcal{B}[\theta_N^{(1)},\theta_N^{(1)}],
\end{align}
and observe that the $\|\theta_N^{(1)}\|_{\dB_{p,q}^{\frac{2}{p}-1}}$ vanishes as $N \to \infty$ while the norm of $\theta_N^{(2)}$ is bounded from below by a positive quantity that does not tends to $0$ as $N \to \infty$.
Then, we write the solution $\theta_N=\theta_N^{(1)}+\theta_N^{(2)}+\widetilde{\theta}_N$, where the perturbation $\widetilde{\theta}_N$ should solve 
\begin{align}\label{eq:pQG}
    \begin{split}
    \widetilde{\theta}_N
    ={}&
    \mathcal{B}[\theta_N^{(1)}, \theta_N^{(2)}]
    +
    \mathcal{B}[\theta_N^{(2)}, \theta_N^{(1)}]
    +
    \mathcal{B}[\theta_N^{(2)}, \theta_N^{(2)}]
    \\
    &
    +
    \mathcal{B}[\theta_N^{(1)}, \widetilde{\theta}_N]
    +
    \mathcal{B}[\widetilde{\theta}_N, \theta_N^{(1)}]
    +
    \mathcal{B}[\theta_N^{(2)}, \widetilde{\theta}_N]
    +
    \mathcal{B}[\widetilde{\theta}_N, \theta_N^{(2)}]
    +
    \mathcal{B}[\widetilde{\theta}_N, \widetilde{\theta}_N].
    \end{split}
\end{align}
Solving \eqref{eq:pQG} in $\dB_{4,2}^{-\frac{1}{2}}(\mathbb{R}^2)$ similarly as in Theorem \ref{thm:WP}, we find that $\widetilde{\theta}_N$ is relatively smaller than $\theta_N^{(2)}$.
Then, combining the estimates for $\theta_N^{(1)}$, $\theta_N^{(2)}$, and $\widetilde{\theta}_N$, we will complete the proof.

\noindent
{\it Step 1. The case of $p >4$ and $1 \leq q \leq \infty$.}
Let $\chi \in \mathscr(\mathbb{R}^2)$ be a radial symmetric real valued function satisfying 
\begin{align}
    \widehat{\chi}(\xi)
    =
    \begin{cases}
    1 & (|\xi| \leq 1),\\
    0 & (|\xi| \geq 2).
    \end{cases}
\end{align}
We define
\begin{align}
    f_N(x)
    :=
    \delta 
    2^{\frac{5}{2}N}
    \chi(x)
    \cos (2^Nx_1), \qquad N \geq 100.
\end{align}
Then, it follows from 
\begin{align}
    \widehat{f_N}(\xi)
    =
    \frac{\delta 2^{\frac{5}{2}N}}{2}
    \sp{
    \widehat{\chi}(\xi+2^Ne_1)
    +
    \widehat{\chi}(\xi-2^Ne_1)
    }
\end{align}
that $\supp \widehat{f_N} \subset \{ \xi \in \mathbb{R}^2\ ;\ 2^N - 2 \leq |\xi| \leq 2^N + 2 \}$, which implies 
\begin{align}
    &
    \n{f_N}_{\dB_{p,q}^{\frac{2}{p}-3}}
    \leq
    C_1
    \delta 
    2^{N(\frac{2}{p}-\frac{1}{2})},\\
    &
    \n{\theta_N^{(1)}}_{\dB_{p,q}^{\frac{2}{p}-1}}
    \leq 
    C
    \n{f_N}_{\dB_{p,q}^{\frac{2}{p}-3}}
    \leq
    C_1
    \delta 
    2^{N(\frac{2}{p}-\frac{1}{2})}
\end{align}
for some positive constant $C_1=C_1(\|\phi_0\|_{L^1},\| \chi \|_{L^p})$.
As $p>4$, we see that 
\begin{align}
    \lim_{N\to \infty}
    \n{f_N}_{\dB_{p,q}^{\frac{2}{p}-3}}
    =0.
\end{align}
Next, we consider the estimates of the second iteration $\theta_N^{(2)}$.
In what follows, we assume $|\xi| \leq 1$.
We see that 
\begin{align}
    &
    \widehat{\theta_N^{(2)}}(\xi)
    =
    \frac{\delta^22^{5N}}{2}
    \frac{\xi_1\xi_2}{|\xi|^2}
    \int_{\mathbb{R}^2}
    \frac{\eta_1^2}{|\xi-\eta|+|\eta|}
    \frac{\widehat{\chi}(\xi-\eta+2^Ne_1)}{|\xi-\eta|^{3}}
    \frac{\widehat{\chi}(\eta-2^Ne_1)}{|\eta|^3}
    d\eta\\
    &\qquad
    +
    \frac{\delta^22^{5N}}{2|\xi|^2}
    \int_{\mathbb{R}^2}
    \frac{(\xi_2^2-\xi_1^2)\eta_1\eta_2-\xi_1\xi_2\eta_2^2}{|\xi-\eta|+|\eta|}
    \frac{\widehat{\chi}(\xi-\eta+2^Ne_1)}{|\xi-\eta|^{3}}
    \frac{\widehat{\chi}(\eta-2^Ne_1)}{|\eta|^3}
    d\eta\\
    &\qquad
    -
    \frac{\delta^22^{5N}}{4}
    \int_{\mathbb{R}^2}
    \frac{ \xi \cdot \eta^{\perp}  }{|\xi-\eta|+|\eta|}
    \frac{\widehat{\chi}(\xi-\eta+2^Ne_1)}{|\xi-\eta|^{3}}
    \frac{\widehat{\chi}(\eta-2^Ne_1)}{|\eta|^3}
    d\eta\\
    &\quad=:
    \widehat{I_N^{(1)}}(\xi)
    +
    \widehat{I_N^{(2)}}(\xi)
    +
    \widehat{I_N^{(3)}}(\xi).
\end{align}
We consider the estimate of $\widehat{I_N^{(1)}}(\xi)$.
We introduce $\psi _n \in \mathscr{S}(\mathbb{R}^2)$ defined by 
\begin{align}
\widehat{\psi_j} (\xi) &= 1- \sum_{k = j-10}^{\infty} \widehat{\phi_k} \big(\xi- 2^ja\big), \quad
a:=
\begin{pmatrix}
    1/\sqrt{2} \\ 1/\sqrt{2}
\end{pmatrix},
\end{align}
which enable us to focus on the frequency in the region $\{ \xi=(\xi_1,\xi_2)\ ;\  \xi_1, \xi _2 > 0\}$.
Note that $\psi_j$ has the following properties:
\begin{align}
    &
    \psi_j(x)=2^{2j}\psi_0(2^jx), \\
    &
    \supp \widehat{\psi_j} \subset \{ \xi \in \mathbb{R}^2\ ;\ |\xi-2^ja| \leq 2^{j-10} \},\\
    &
    \psi_j = \phi_j *\psi_j
\end{align}
for $j\in \mathbb{Z}$.
In particular, it follows from the last property that 
\begin{align}
    \n{\psi_j *f}_{L^p} \leq C \n{\phi_j *f}_{L^p}
\end{align} 
for all $1 \leq p \leq \infty$, $j \in \mathbb{Z}$, and distributions $f$ with $\Delta_j f \in L^p(\mathbb{R}^2)$.
For $\eta \in \supp \sp{\widehat{\chi}(\xi-\cdot+2^Ne_1)\widehat{\chi}(\cdot-2^Ne_1)}$,
we see that $|\xi| \ll |\xi-\eta| \sim |\eta| \sim \eta_1 \sim 2^N$ and $|\eta_2| \lesssim 1$.  
Thus, we have
\begin{align}
    \widehat{\psi_j}(\xi)
    \widehat{I_N^{(1)}}(\xi) 
    \sim{}& 
    \delta^2
    \widehat{\psi_j}(\xi)
    \frac{\xi_1\xi_2}{|\xi|^2}
    \int_{\mathbb{R}^2}
    \widehat{\chi}(\xi-\eta+2^Ne_1)
    \widehat{\chi}(\eta-2^Ne_1)
    d\eta\\
    ={}& 
    \delta^2
    \widehat{\psi_j}(\xi)
    \frac{\xi_1\xi_2}{|\xi|^2}
    \int_{\mathbb{R}^2}
    {\widehat{\chi}(\xi-\eta')}
    {\widehat{\chi}(\eta')}d\eta'\\
    \sim{}& 
    \delta^2
    \widehat{\psi_j}(\xi).
\end{align}
For the estimate of $\widehat{I_N^{(2)}}(\xi)$, it holds 
\begin{align}
    \left| \widehat{I_N^{(2)}}(\xi) \right|
    \leq{}&
    C\delta^22^{5N}
    \int_{\mathbb{R}^2}
    \frac{\eta_1|\eta_2|}{|\xi-\eta|+|\eta|}
    \frac{\widehat{\chi}(\xi-\eta+2^Ne_1)}{|\xi-\eta|^{3}}
    \frac{\widehat{\chi}(\eta-2^Ne_1)}{|\eta|^3}
    d\eta\\
    \leq{}&
    C\delta^22^{-N}
    \int_{\mathbb{R}^2}
    \widehat{\chi}(\xi-\eta+2^Ne_1)
    \widehat{\chi}(\eta-2^Ne_1)
    d\eta\\
    \leq{}&
    C\delta^22^{-N}.
\end{align}
For the estimate of $\widehat{I_N^{(3)}}(\xi)$, we have
\begin{align}
    \left|\widehat{I_N^{(3)}}(\xi)\right|
    \leq{}&
    \delta^22^{5N}|\xi|
    \int_{\mathbb{R}^2}
    \frac{|\eta|}{|\xi-\eta|+|\eta|}
    \frac{\widehat{\chi}(\xi-\eta+2^Ne_1)}{|\xi-\eta|^{3}}
    \frac{\widehat{\chi}(\eta-2^Ne_1)}{|\eta|^3}
    d\eta\\
    \leq{}&
    C\delta^22^{-N}|\xi|
    \int_{\mathbb{R}^2}
    \widehat{\chi}(\xi-\eta+2^Ne_1)
    \widehat{\chi}(\eta-2^Ne_1)
    d\eta\\
    \leq{}&
    C\delta^22^{-N}.
\end{align}
Hence, we obtain for each $j \leq -1$ that
\begin{align}
    &2^{-j}
    \n{\phi_j * \theta_N^{(2)}}_{L^{\infty}}
    \geq{}
    c
    2^{-j}
    \n{\psi_j * I_N^{(1)}}_{L^{\infty}}
    -
    \sum_{m=2}^3
    2^{-j}
    \n{\phi_j * I_N^{(m)}}_{L^{\infty}}
    \\
    &\quad\geq{}
    c
    \n{e^{-2^{2j}|\cdot|^2}\sp{\psi_j * I_N^{(1)}}}_{L^{2}}
    -
    C
    \sum_{m=2}^3
    \n{\phi_j * I_N^{(m)}}_{L^{2}}
    \\
    &\quad={}
    c
    \n{\frac{1}{4\pi2^{2j}}\int_{\mathbb{R}^2}e^{-\frac{|\xi-\xi'|^2}{2^{2j+2}}}\widehat{\psi_0}(2^{-j}\xi') \widehat{I_N^{(1)}}(\xi')d\xi'}_{L^2}
    -
    C
    \sum_{m=2}^3
    \n{\widehat{\phi_j}(\xi) \widehat{I_N^{(m)}}(\xi)}_{L^{2}}\\
    &\quad\geq{}
    c\delta^2
    2^{-2j}
    \n{\int_{\mathbb{R}^2}e^{-\frac{|\xi-\xi'|^2}{2^{2j+2}}}\widehat{\psi_0}(2^{-j}\xi')d\xi'}_{L^2}-C\delta^22^{-N}\n{\widehat{\phi_0}(2^{-j}\xi)}_{L^2}\\
    &\quad={}
    c\delta^2
    \n{\int_{\mathbb{R}^2}e^{-\frac{|2^{-j}\xi-\xi'|^2}{4}}\widehat{\psi_0}(\xi')d\xi'}_{L^2}-C\delta^22^j2^{-N}
    \\
    &\quad={}
    c\delta^22^j
    \n{\int_{\mathbb{R}^2}e^{-\frac{|\xi-\xi'|^2}{4}}\widehat{\psi_0}(\xi')d\xi'}_{L^2}-C\delta^22^j2^{-N}\\
    &\quad={}
    c\delta^22^j-C\delta^22^j2^{-N}.
\end{align}

If $\theta_N$ is a solution to \eqref{eq:QG} with the external force $f_N$, then the higher iteration part $\widetilde{\theta}_N:=\theta_N - \theta_N^{(1)} - \theta_N^{(2)}$ should solve \eqref{eq:pQG}.
From Lemma \ref{lemm:nonlin-1} it follows that
\begin{align}
    \n{\mathcal{B}[\theta_N^{(1)}, \theta_N^{(2)}]
    +
    \mathcal{B}[\theta_N^{(2)}, \theta_N^{(1)}]}_{\dB_{4,2}^{-\frac{1}{2}}}
    &\leq
    C
    \n{\theta_N^{(1)}}_{\dB_{4,2}^{-\frac{1}{2}}}
    \n{\theta_N^{(2)}}_{\dB_{4,2}^{-\frac{1}{2}}}
    \leq
    C\delta^3,\\
    \n{\mathcal{B}[\theta_N^{(2)},\theta_N^{(2)}]}_{\dB_{4,2}^{-\frac{1}{2}}}
    &\leq
    C
    \n{\theta_N^{(2)}}_{\dB_{4,2}^{-\frac{1}{2}}}^2
    \leq
    C\delta^4,\\
    \n{\mathcal{B}[\theta_N^{(1)}, \theta]
    +
    \mathcal{B}[\theta, \theta_N^{(1)}]}_{\dB_{4,2}^{-\frac{1}{2}}}
    &\leq
    C
    \n{\theta_N^{(1)}}_{\dB_{4,2}^{-\frac{1}{2}}}
    \n{\theta}_{\dB_{4,2}^{-\frac{1}{2}}}
    \leq
    C\delta
    \n{\theta}_{\dB_{4,2}^{-\frac{1}{2}}},\\
    \n{\mathcal{B}[\theta_N^{(2)}, {\theta}]
    +
    \mathcal{B}[\theta, \theta_N^{(2)}]}_{\dB_{4,2}^{-\frac{1}{2}}}
    &\leq
    C
    \n{\theta_N^{(2)}}_{\dB_{4,2}^{-\frac{1}{2}}}
    \n{\theta}_{\dB_{4,2}^{-\frac{1}{2}}}\\
    &\leq
    C
    \n{\theta_N^{(1)}}_{\dB_{4,2}^{-\frac{1}{2}}}^2
    \n{{\theta}}_{\dB_{4,2}^{-\frac{1}{2}}}\\
    &\leq
    C\delta^2
    \n{{\theta}}_{\dB_{4,2}^{-\frac{1}{2}}},\\
    \n{\mathcal{B}[{\theta}, \theta']}_{\dB_{4,2}^{-\frac{1}{2}}}
    &\leq
    C
    \n{\theta}_{\dB_{4,2}^{-\frac{1}{2}}}
    \n{\theta'}_{\dB_{4,2}^{-\frac{1}{2}}}
\end{align}
for all $\theta,\theta' \in \dB_{4,2}^{-\frac{1}{2}}(\mathbb{R}^2)$.
Then, it follows from the similar contraction mapping argument via the suitable small $\delta$ as in the proof of Theorem \ref{thm:WP} that there exists a unique solution $\widetilde{\theta}_N \in \dB_{4,2}^{-\frac{1}{2}}(\mathbb{R}^2)$ to \eqref{eq:pQG} satisfying 
\begin{align}
    \n{\widetilde{\theta}_N}_{\dB_{4,2}^{-\frac{1}{2}}}
    \leq
    C\delta^3
\end{align}
with some positive constant $C$ independent of $N$.
Hence, $\theta_N:=\theta_N^{(1)} + \theta_N^{(2)} + \widetilde{\theta}_N$ solves \eqref{eq:QG} with the force $f_N$ and satisfies 
\begin{align}
    \n{\theta_N}_{\dB_{\infty,\infty}^{-1}}
    \geq{}&
    \sup_{j \leq -1}
    2^{-j}
    \n{\phi_j*\theta_N^{(2)}}_{L^{\infty}}
    -
    \n{\theta_N^{(1)}}_{\dB_{\infty,\infty}^{-1}}
    -
    C
    \n{\widetilde{\theta}_N}_{\dB_{4,2}^{-\frac{1}{2}}}\\
    \geq{}&
    c\delta^2 - C2^{-N} - C\delta 2^{-\frac{N}{2}} - C\delta^3.
\end{align}
choosing $\delta$ sufficiently small, we obtain $\liminf_{N\to \infty}\n{\theta_N}_{\dB_{\infty,\infty}^{-1}} > 0$.

\noindent
{\it Step 2. The case of $p =4$ and $2 < q \leq \infty$.}
Let 
\begin{align}
    f_N(x)
    :=
    \frac{\delta}{\sqrt{\log N}}
    \sum_{n=10}^N  
    \frac{2^{\frac{5}{2}n^2}}{\sqrt{n}}
    \chi(x)
    \cos (2^{n^2}x_1),\qquad N \geq 100.
\end{align}
Then, we see that 
\begin{align}
    \n{f_N}_{\dB_{4,q}^{-\frac{5}{2}}}
    \leq 
    \frac{C\delta}{\sqrt{\log N}}
    \Mp{\sum_{n=10}^N
    \frac{1}{n^{\frac{q}{2}}}}^{\frac{1}{q}}
    \leq
    \begin{cases}
    \dfrac{C\delta}{\sqrt{\log N}} & (2<q\leq \infty),\\
    C\delta & (q=2)
    \end{cases}
\end{align}
and
\begin{align}
    \n{\theta_N^{(1)}}_{\dB_{4,q}^{-\frac{1}{2}}}
    \leq
    C 
    \n{f_N}_{\dB_{4,q}^{-\frac{5}{2}}}
    \leq
    \begin{cases}
    \dfrac{C\delta}{\sqrt{\log N}} & (2<q\leq \infty),\\
    C\delta & (q=2).
    \end{cases}
\end{align}
Similarly to the above calculations as in Step 1, we have
\begin{align}
    \widehat{\theta_N^{(2)}}(\xi)
    &={}
    \sum_{n=M}^N
    \frac{1}{n{\log N}}
    \lp{\widehat{I_{n^2}^{(1)}}(\xi)+\widehat{I_{n^2}^{(2)}}(\xi)+\widehat{I_{n^2}^{(3)}}(\xi)}\\
    &=:{}
    \widehat{J_{N}^{(1)}}(\xi)+\widehat{J_{N}^{(2)}}(\xi)+\widehat{J_{N}^{(3)}}(\xi)
\end{align}
for all $\xi \in \mathbb{R}^2$ with $|\xi| \leq 1$.
Then, we see that 
\begin{align}
    \widehat{\psi_j}(\xi)
    \widehat{J_N^{(1)}}(\xi)
    =
    \frac{1}{\log N}
    \sum_{n=10}^N
    \frac{1}{n}
    \widehat{\psi_j}(\xi)
    \widehat{I_{n^2}^{(1)}}(\xi)
    \sim
    -
    \delta^2
    \widehat{\psi_j}(\xi)
\end{align}
and 
\begin{align}
    \sum_{m=2}^3
    \left| \widehat{J_N^{(m)}}(\xi) \right|
    \leq{}&
    \sum_{m=2}^3
    \sum_{n=10}^N
    \frac{1}{n\log N}
    \left| \widehat{I_N^{(m)}}(\xi) \right|\\
    \leq{}&
    C
    \frac{\delta^2}{\log N}
    \sum_{n=10}^N
    \frac{2^{-n}}{n}
    \\
    \leq{}&
    C
    \frac{\delta^2}{\log N}.
\end{align}
Then, from the similar calculation as above, it holds that
\begin{align}
    &2^{-j}
    \n{\phi_j * \theta_N^{(2)}}_{L^{\infty}}
    \geq{}
    c\n{e^{-2^{2j}|\cdot|^2}\sp{\psi_j * J_N^{(1)}}}_{L^2}
    -
    C
    \sum_{m=2}^3
    \n{\phi_j * J_N^{(m)}}_{L^2}
    \\
    &\quad
    \geq{}
    c\delta^2
    \n{\frac{1}{4\pi2^{2j}}\int_{\mathbb{R}^2}e^{-\frac{|\xi-\xi'|^2}{2^{2j+2}}}\widehat{\psi_0}(2^{-j}\xi') \widehat{J_N^{(1)}}(\xi')d\xi'}_{L^2}
    -
    C
    \sum_{m=2}^3
    \n{\widehat{\phi_j}(\xi) \widehat{J_N^{(m)}}(\xi) }_{L^2}
    \\
    &\quad 
    \geq{}
    c\delta^22^j-C\frac{\delta^2}{\log N}2^j
\end{align}
for all $j \geq -1$.
The rest of the proof is completely same as in Step 1.

\noindent
{\it Step 3. The case of $p =4$ and $1 \leq q < 2$.}
Let $N \geq 200$ and let us define 
\[
f_N(x) = \dfrac{\delta 2^{\frac{5}{2}N^2}}{N^\frac{1}{4}\log N} F_N(x)\cos \sp{2^{N^2} x_1},
\]
where we have set
\begin{align}
    F_N(x)
    :=
    \sum _{k=10}^{N-50}
    2^{-\frac{3}{2}k^2} \phi_{k^2} (x - Rk^2e_1).
\end{align}
Here, we choose $R>0$ so large that 
\begin{align}
    \n{F_N}_{L^4} \leq CN^{\frac{1}{4}}
\end{align}
holds with some constant $C$ independent of $N$, which is proved by the similar calculation as in \cites{CN,IO}.
We notice that 
\begin{align}
    \supp \widehat{f_N} &\subset \Mp{ \xi \in \mathbb{R}^2 \ ;\ 2^{N^2-1} \leq |\xi| \leq 2^{N^2+1} },\\
    \supp \widehat{F_N} &\subset \Mp{ \xi \in \mathbb{R}^2 \ ;\ 2^{99} \leq |\xi| \leq 2^{(N-50)^2+1} }.
\end{align}
For the estimate of $f_N$, it holds  
\begin{align}
    \| f_N \|_{\dot B^{-\frac{5}{2}}_{4,\sigma}} 
    \leq{}&
    \left\{ 
    \sum_{j=N^2-1}^{N^2+1}
    \left(
    2^{-\frac{5}{2}j}
    \n{\phi_j * f_N}_{L^4}
    \right)^\sigma
    \right\}^{\frac{1}{\sigma}}
    \\
    \leq{}&  
    \dfrac{C\delta}{N^\frac{1}{4}\log N}
    \left\| 
        F_N
    \right\|_{L^4}
    \leq 
    \frac{C\delta}{\log N}
\end{align}
for all $1 \leq \sigma \leq \infty$.
Similarly as in the previous steps, we focus on the first iteration $\theta_N^{(1)}$ and the second iteration $\theta_N^{(2)}$.
For the estimate of $\theta_N^{(1)}$, we have 
\begin{align}\label{est:1}
    \n{\theta_N^{(1)}}_{\dB_{4,2}^{-\frac{1}{2}}}
    \leq 
    C
    \n{f_N}_{\dB_{4,2}^{-\frac{5}{2}}}
    \leq 
    C \frac{\delta}{\log N}.
\end{align}
We next consider the estimate of the norm of $\theta_N^{(2)}$ in $\dot B^{-\frac{1}{2}}_{4,q}(\mathbb{R}^2)$.  
Similarly to Step 1, we decompose $\theta_N^{(2)}$ as 
{\small\begin{align}
    &\widehat{\theta_N^{(2)}}(\xi)
    =
    \frac{\delta^22^{5N^2}}{2N^{\frac{1}{2}}(\log N)^2}
    \frac{\xi_1\xi_2}{|\xi|^2}
    \int_{\mathbb{R}^2}
    \frac{\eta_1^2}{|\xi-\eta|+|\eta|}
    \frac{\widehat{F_N}(\xi-\eta+2^{N^2}e_1)}{|\xi-\eta|^{3}}
    \frac{\widehat{F_N}(\eta-2^{N^2}e_1)}{|\eta|^3}
    d\eta\\
    &
    +
    \frac{\delta^22^{5N^2}}{2N^{\frac{1}{2}}(\log N)^2}
    \int_{\mathbb{R}^2}
    \frac{(\xi_2^2-\xi_1^2)\eta_1\eta_2-\xi_1\xi_2\eta_2^2}{|\xi|^2(|\xi-\eta|+|\eta|)}
    \frac{\widehat{F_N}(\xi-\eta+2^{N^2}e_1)}{|\xi-\eta|^{3}}
    \frac{\widehat{F_N}(\eta-2^{N^2}e_1)}{|\eta|^3}
    d\eta\\
    & 
    -
    \frac{\delta^22^{5N^2}}{4N^{\frac{1}{2}}(\log N)^2}
    \int_{\mathbb{R}^2}
    \frac{ \xi \cdot \eta^{\perp}  }{|\xi-\eta|+|\eta|}
    \frac{\widehat{F_N}(\xi-\eta+2^{N^2}e_1)}{|\xi-\eta|^{3}}
    \frac{\widehat{F_N}(\eta-2^{N^2}e_1)}{|\eta|^3}
    d\eta\\
    &=:
    \widehat{K_N^{(1)}}(\xi)
    +
    \widehat{K_N^{(2)}}(\xi)
    +
    \widehat{K_N^{(3)}}(\xi)
\end{align}}
for all $\xi$ with $|\xi| \leq 2^{(N-100)^2}$.
We first consider the estimate of $K_N^{(1)}$.
If
\begin{align}
    \supp \sp{\widehat{\psi_{n^2}}(\xi)\widehat{\phi_{k^2}}(\cdot)\widehat{\phi_{\ell^2}}(\xi-\cdot)}
    \neq
    \varnothing
\end{align}
then it should hold
$(k,\ell) \in ([10,n-1] \times \{n\}) \cup (\{n\} \times [10,n-1]) \cup \{ (k,k)\ ;\ k\geq n \}$,
which implies
\begin{align}
    &
    \widehat{\psi_{n^2}}(\xi)
    \widehat{F_N}(\xi-\eta+2^{N^2}e_1)
    \widehat{F_N}(\eta-2^{N^2}e_1)
    \\
    &\quad 
    ={}
    \widehat{\psi_{n^2}}(\xi)
    \sum_{k,\ell=10}^{N-50}
    \sp{\begin{aligned}
    &
    e^{-iRk^2(\xi_1-\eta_1)}
    e^{-iR\ell^2\eta_1}
    2^{-\frac{3}{2}k^2}
    2^{-\frac{3}{2}n^2}\\
    &
    \times 
    \widehat{\phi_{k^2}}(\xi-\eta+2^{N^2}e_1)
    \widehat{\phi_{n^2}}(\eta-2^{N^2}e_1)
    \end{aligned}}
    \\
    &\quad 
    ={}
    \widehat{\psi_{n^2}}(\xi)
    \sum_{k=n}^{N-50}
    e^{-iRk^2\xi_1}
    2^{-3k^2}
    \widehat{\phi_{k^2}}(\xi-\eta+2^{N^2}e_1)
    \widehat{\phi_{k^2}}(\eta-2^{N^2}e_1)\\
    &\qquad
    +
    \widehat{\psi_{n^2}}(\xi)
    \sum_{\ell=10}^{n-1}
    \sp{\begin{aligned}
    &
    e^{-iRn^2(\xi_1-\eta_1)}
    e^{-iR\ell^2\eta_1}
    2^{-\frac{3}{2}n^2}
    2^{-\frac{3}{2}\ell^2}\\
    &
    \times 
    \widehat{\phi_{n^2}}(\xi-\eta+2^{N^2}e_1)
    \widehat{\phi_{\ell^2}}(\eta-2^{N^2}e_1)
    \end{aligned}}
    \\
    &\qquad
    +
    \widehat{\psi_{n^2}}(\xi)
    \sum_{k=10}^{n-1}
    \sp{\begin{aligned}
    &
    e^{-iRk^2(\xi_1-\eta_1)}
    e^{-iRn^2\eta_1}
    2^{-\frac{3}{2}k^2}
    2^{-\frac{3}{2}n^2}\\
    &
    \times 
    \widehat{\phi_{k^2}}(\xi-\eta+2^{N^2}e_1)
    \widehat{\phi_{n^2}}(\eta-2^{N^2}e_1)
    \end{aligned}}.
\end{align}
Thus, we see that 
\begin{align}
    &
    \psi_n * K_N^{(1)}(x)\\
    &\quad
    =
    -
    \frac{\delta^22^{5N^2}}{2N^{\frac{1}{2}}(\log N)^2}
    2^{-3n^2} 
    \frac{1}{(2\pi)^2}
    \int_{\mathbb{R}^2}
    e^{i\xi \cdot (x-Rn^2e_1)}
    G_{n,n,N}(\xi)
    d\xi\\
    &
    \qquad
    -
    \frac{\delta^22^{5N^2}}{2N^{\frac{1}{2}}(\log N)^2}
    \sum_{k=n+1}^{N-50}
    2^{-3k^2} 
    \frac{1}{(2\pi)^2}
    \int_{\mathbb{R}^2}
    e^{i\xi \cdot (x-Rk^2e_1)}
    G_{n,k,N}(\xi)
    d\xi\\
    &
    \qquad
    -
    \frac{\delta^22^{5N^2}}{2N^{\frac{1}{2}}(\log N)^2}
    \sum_{\ell=10}^{n-1}
    2^{-\frac{3}{2}n^2}
    2^{-\frac{3}{2}\ell^2}
    \frac{1}{(2\pi)^2}
    \int_{\mathbb{R}^2}
    e^{i\xi \cdot (x-Rk^2e_1)}
    H_{n,n,\ell,N}(\xi)
    d\xi\\
    &
    \qquad
    -
    \frac{\delta^22^{5N^2}}{2N^{\frac{1}{2}}(\log N)^2}
    \sum_{k=10}^{n-1}
    2^{-\frac{3}{2}k^2}
    2^{-\frac{3}{2}\ell^2}
    \frac{1}{(2\pi)^2}
    \int_{\mathbb{R}^2}
    e^{i\xi \cdot (x-Rk^2e_1)}
    H_{n,k,n,N}(\xi)
    d\xi\\
    &\quad 
    =:K_{N,n}^{(1,1)}(x)+K_{N,n}^{(1,2)}(x)+K_{N,n}^{(1,3)}(x)+K_{N,n}^{(1,4)}(x),
\end{align}
where we have set 
\begin{align}
    G_{n,k,N}(\xi)
    :={}&
    \frac{\xi_1\xi_2}{|\xi|^2}
    \widehat{\psi_{n^2}}(\xi)
    \int_{\mathbb{R}^2}
    \frac{\eta_1^2}{|\xi-\eta|+|\eta|}
    \frac{\widehat{\phi_{k^2}}(\xi-\eta+2^{N^2}e_1)}{|\xi-\eta|^3}
    \frac{\widehat{\phi_{k^2}}(\eta-2^{N^2}e_1)}{|\eta|^3}
    d\eta,\\
    H_{n,k,\ell,N}(\xi)
    :={}&
    \frac{\xi_1\xi_2}{|\xi|^2}
    \widehat{\psi_{n^2}}(\xi)
    \int_{\mathbb{R}^2}
    e^{iR(k^2-\ell^2)\eta_1}
    \frac{\eta_1^2}{|\xi-\eta|+|\eta|}\\
    &\qquad \qquad \qquad \times
    \frac{\widehat{\phi_{k^2}}(\xi-\eta+2^{N^2}e_1)}{|\xi-\eta|^3}
    \frac{\widehat{\phi_{\ell^2}}(\eta-2^{N^2}e_1)}{|\eta|^3}
    d\eta.
\end{align}
Let us consider the estimate of $K_{N}^{(1,1)}(x)$.
We remark that $G_{n,k,N}(\xi) \geq 0$ holds.
Let  
\begin{align}
    A_{n}:= \Mp{ x \in \mathbb{R}^2 \ ;\ |x-Rn^2e_1| \leq \varepsilon2^{-n^2} },
\end{align}
where $\varepsilon$ is a positive constant to be determined later.
Then, we have
\begin{align}
    \left| e^{i\xi \cdot (x-Rn^2e_1)} - 1 \right| \leq |\xi||x-Rn^2e_1| \leq 2\varepsilon, \qquad
    x \in A_{n}, \quad \xi \in \supp \widehat{\psi_n}.
\end{align}
and $|A_{n}| \leq C\varepsilon^22^{-2n^2}$.
Thus, it holds 
\begin{align}
    &
    \n{K_{N,n}^{(1,1)}}_{L^4(A_{n})}\\
    &\quad
    \geq{}
    c\frac{\delta^22^{5N^2}}{N^{\frac{1}{2}}(\log N)^2}
    2^{-3n^2}
    |A_{n}|^{\frac{1}{4}}
    \int_{\mathbb{R}^2}
    G_{n,n,N}(\xi)
    d\xi\\
    &\qquad
    -
    C
    \frac{\delta^22^{5N^2}}{N^{\frac{1}{2}}(\log N)^2}
    2^{-3n^2} 
    \n{\int_{\mathbb{R}^2}
    \left| e^{i\xi \cdot (x-Rn^2e_1)}-1 \right|
    G_{n,n,N}(\xi)
    d\xi}_{L^4(A_{n})} 
    \\
    &\quad
    \geq{}
    c
    \frac{\delta^22^{5N^2}}{N^{\frac{1}{2}}(\log N)^2}
    2^{-3n^2}
    \sp{\varepsilon^2 2^{-2n^2}}^{\frac{1}{4}}
    \int_{\mathbb{R}^2}
    G_{n,n,N}(\xi)
    d\xi\\
    &\qquad
    -
    C
    \frac{\delta^22^{5N^2}}{N^{\frac{1}{2}}(\log N)^2}
    2^{-3n^2} 
    \int_{\mathbb{R}^2}
    \varepsilon 2^{-n^2}
    G_{n,n,N}(\xi)
    d\xi \\
    &\quad
    \geq{}
    \left(c_*\varepsilon^{\frac{1}{2}}-C_*\varepsilon\right)
    \frac{\delta^22^{5N^2}}{N^{\frac{1}{2}}(\log N)^2}
    2^{-3n^2}
    2^{-\frac{1}{2}n^2}
    \int_{\mathbb{R}^2}
    G_{n,n,N}(\xi)
    d\xi
\end{align}
with some positive constants $c_*$ and $C_*$.
We choose $\varepsilon$ so small that $c_*\varepsilon^{\frac{1}{2}}-C_*\varepsilon > 0$.
Here, using
\begin{align}
    G_{n,k,N}(\xi)
    \sim{}&
    \widehat{\psi_{n^2}}(\xi)
    \int_{\mathbb{R}^2}
    \frac{2^{2N^2}}{2^{N^2}}
    \frac{\widehat{\phi_{k^2}}(\xi-\eta+2^{N^2}e_1)}{2^{3N^2}}
    \frac{\widehat{\phi_{k^2}}(\eta-2^{N^2}e_1)}{2^{3N^2}}
    d\eta \\
    \sim{}&
    2^{-5N^2}
    \widehat{\psi_{n^2}}(\xi)
    \int_{\mathbb{R}^2}
    \widehat{\phi_{k^2}}(\xi-\eta)\widehat{\phi_{k^2}}(\eta)
    d\eta\\ 
    \sim{}&
    2^{-5N^2}
    \widehat{\psi_{n^2}}(\xi)
    \sp{\widehat{\phi_{k^2}} * \widehat{\phi_{k^2}}}(\xi),
\end{align}
we see that
\begin{align}
    \int_{\mathbb{R}^2}
    G_{n,k,N}(\xi)
    d\xi 
    \sim {}&
    2^{-5N^2}
    \int_{\mathbb{R}^2}
    \widehat{\psi_{n^2}}(\xi)
    \sp{\widehat{\phi_{k^2}} * \widehat{\phi_{k^2}}}(\xi)
    d\xi\\
    \sim {}&
    2^{-5N^2}2^{2k^2}
    \int_{\mathbb{R}^2}
    \widehat{\psi_0}(2^{-n^2}\xi)
    \sp{\widehat{\phi_0} * \widehat{\phi_0}}(2^{-k^2}\xi)
    d\xi\\
    \sim {}&
    2^{-5N^2}2^{2n^2}2^{2k^2}
\end{align}
for $k\geq n \geq 100$.
Thus, we have 
\begin{align}\label{K11}
    \n{K_{N,n}^{(1,1)}}_{L^4(A_{n})}
    \geq
    c\varepsilon^{\frac{1}{2}}
    2^{\frac{1}{2}n^2}
    \frac{\delta^2}{N^{\frac{1}{2}}(\log N)^2}.
\end{align}
For the estimate of $K_{N,n}^{(1,2)}$,
it holds
\begin{align}\label{K12}
    \begin{split}
    \n{K_{N,n}^{(1,2)}}_{L^4(A_n)}
    \leq{}&
    C
    \frac{\delta^22^{5N^2}}{2N^{\frac{1}{2}}(\log N)^2}
    \sum_{k=n+1}^{N-50}
    2^{-3k^2} 
    \sp{\varepsilon^22^{-2n^2}}^{\frac{1}{4}}
    \int_{\mathbb{R}^2}
    G_{n,k,N}(\xi)
    d\xi\\
    \leq{}&
    C2^{\frac{3}{2}n^2}
    \frac{\delta^2}{2N^{\frac{1}{2}}(\log N)^2}
    \sum_{k=n+1}^{N-50}
    2^{-k^2} 
    \\
    \leq{}&
    C2^{\frac{3}{2}n^2}
    2^{-(n+1)^2}
    \frac{\delta^2}{2N^{\frac{1}{2}}(\log N)^2}
    =
    C2^{\frac{1}{2}n^2}
    2^{-2n-1}
    \frac{\delta^2}{2N^{\frac{1}{2}}(\log N)^2}.
    \end{split}
\end{align}
For the estimate of $K_{N,n}^{(1,3)}$, it holds
\begin{align}
    \n{K_{N,n}^{(1,3)}}_{L^4(A_n)}
    &
    \leq 
    C
    \frac{\delta^22^{5N^2}}{N^{\frac{1}{2}}(\log N)^2}
    \sum_{\ell=10}^{n-1}
    2^{-\frac{3}{2}(k^2+\ell^2)}
    \n{
    \int_{\mathbb{R}^2}
    e^{i(x-Rk^2e_1)\cdot \xi}
    H_{n,n,\ell,N}(\xi)
    d\xi
    }_{L^4(A_{n})}\\
    &
    \leq{}
    C
    \frac{\delta^22^{5N^2}}{N^{\frac{1}{2}}(\log N)^2}
    \sum_{\ell=10}^{n-1}
    2^{-\frac{3}{2}(n^2+\ell^2)}
    \sp{\varepsilon^2 2^{-2n^2}}^{\frac{1}{4}}
    \int_{\mathbb{R}^2}
    \left|H_{n,n,\ell,N}(\xi)\right|
    d\xi\\
    &
    \leq{}
    C
    2^{-\frac{1}{2}n^2}
    \frac{\delta^22^{5N^2}}{N^{\frac{1}{2}}(\log N)^2}
    \sum_{\ell=10}^{n-1}
    2^{-\frac{3}{2}(n^2+\ell^2)}
    \int_{\mathbb{R}^2}
    \left|H_{n,n,\ell,N}(\xi)\right|
    d\xi
\end{align}
Since we see by $\widehat{\psi_{n^2}}(\xi) \leq \widehat{\phi_{n^2}}(\xi)$ that
\begin{align}
    \int_{\mathbb{R}^2}
    \left|H_{n,n,\ell,N}(\xi)\right|
    d\xi 
    \leq{}&
    C2^{-5N^2}
    \int_{\mathbb{R}^2}
    \widehat{\phi_{n^2}}(\xi)
    \int_{\mathbb{R}^2}
    \widehat{\phi_{n^2}}(\xi-\eta+2^{N^2}e_1)
    \widehat{\phi_{\ell^2}}(\eta-2^{N^2}e_1)
    d\eta
    d\xi\\
    ={}&
    C2^{-5N^2}
    \int_{\mathbb{R}^2}\widehat{\phi_{0}}(2^{-n^2}\xi)
    \int_{\mathbb{R}^2}
    \widehat{\phi_0}(2^{-n^2}\xi-2^{-n^2}\eta)
    \widehat{\phi_0}(2^{-\ell^2}\eta)
    d\eta
    d\xi
    \\
    ={}&
    C2^{-5N^2}
    2^{2(n^2+\ell^2)}
    \int_{\mathbb{R}^2}\widehat{\phi_{0}}(\xi)
    \int_{\mathbb{R}^2}
    \widehat{\phi_0}(\xi-2^{\ell^2-n^2}\eta)
    \widehat{\phi_0}(\eta)
    d\eta
    d\xi
    \\
    ={}&
    C2^{-5N^2}
    2^{2(n^2+\ell^2)}
    \int_{\mathbb{R}^2}
    \sp{\widehat{\phi_{0}}*\widehat{\phi_0}}(2^{\ell^2-n^2}\eta)
    \widehat{\phi_0}(\eta)
    d\eta
    \\
    \leq{}&
    C2^{-5N^2}
    2^{2(n^2+\ell^2)}
\end{align}
for $\ell \leq n-1$, 
there holds
\begin{align}\label{K13}
    \begin{split}
    \n{K_{N,n}^{(1,3)}}_{L^4(A_n)}
    \leq{}&
    C
    2^{-\frac{1}{2}n^2}
    \frac{\delta^2}{N^{\frac{1}{2}}(\log N)^2}
    \sum_{\ell=10}^{n-1}
    2^{\frac{1}{2}(n^2+\ell^2)}\\
    \leq{}&
    C
    2^{\frac{1}{2}(n-1)^2}
    \frac{\delta^2}{N^{\frac{1}{2}}(\log N)^2}.
    \end{split}
\end{align}
Similarly, we have
\begin{align}\label{K14}
    \n{K_{N,n}^{(1,4)}}_{L^4(A_n)}
    \leq{}&
    C2^{\frac{1}{2}(n-1)^2}
    \frac{\delta^2}{N^{\frac{1}{2}}(\log N)^2}.
\end{align}
Combining \eqref{K11}, \eqref{K12}, \eqref{K13}, and \eqref{K14}, we have
\begin{align}
    \n{\psi_{n^2}*K_N^{(1)}}_{L^4}
    &\geq 
    \n{\psi_{n^2}*K_N^{(1)}}_{L^4(A_n)}\\
    &\geq 
    \n{K_{N,n}^{(1,1)}}_{L^4(A_n)}
    -
    \n{K_{N,n}^{(1,2)}}_{L^4(A_n)}\\
    &
    \qquad
    -
    \n{K_{N,n}^{(1,3)}}_{L^4(A_n)}
    -
    \n{K_{N,n}^{(1,4)}}_{L^4(A_n)}\\
    &
    \geq
    \sp{c-C2^{-n}}2^{\frac{1}{2}n^2}
    \frac{\delta^2}{N^{\frac{1}{2}}(\log N)^2},
\end{align}
which implies 
\begin{align}
    \left\{ \sum _{n=\frac{N}{2}}^{N-100}
    \left( 2^{-\frac{1}{2}n^2} 
    \left\| \psi_{n^2} * K_N^{(1)} \right\|_{L^4}
    \right)^q \right\}^{\frac{1}{q}}
    \geq
    c
    \frac{\delta^2}{(\log N)^2}
    N^{\frac{1}{q}-\frac{1}{2}}
\end{align}
by choosing $N$ sufficiently large.
Next, we consider the estimates of $K_N^{(m)}$ ($m=2,3$), we have 
\begin{align}
    \sum_{m=2}^3
    \left|\widehat{\psi_{n^2}}(\xi)\widehat{ K_N^{(m)}}(\xi) \right|
    \leq{}&
    C
    \frac{\delta^22^{5N^2}}{N^{\frac{1}{2}}(\log N)^2}
    \widehat{\psi_{n^2}}(\xi)\\
    &\times 
    \int_{\mathbb{R}^2}
    \frac{\widehat{F_N}(\xi-\eta+2^{N^2}e_1)}{|\xi-\eta|^{3}}
    \frac{\widehat{F_N}(\eta-2^{N^2}e_1)}{|\eta|^3}
    d\eta\\
    \leq{}&
    C
    \frac{\delta^22^{-N^2}}{N^{\frac{1}{2}}(\log N)^2}
    \widehat{\psi_{n^2}}(\xi)
    \int_{\mathbb{R}^2}
    \widehat{F_N}(\xi-\eta+2^{N^2}e_1)
    \widehat{F_N}(\eta-2^{N^2}e_1)
    d\eta\\
    ={}&
    C
    \frac{\delta^22^{-N^2}}{N^{\frac{1}{2}}(\log N)^2}
    \widehat{\psi_{n^2}}(\xi)
    \widehat{(F_N)^2}(\xi),
\end{align} 
which yields
\begin{align}
    \Mp{ \sum_{n=\frac{N}{2}}^{N-100} \sp{ 2^{-\frac{1}{2}n^2} \n{\phi_{n^2}*K_N^{(3)}}_{L^4} }^q}^{\frac{1}{q}}
    &\leq
    C
    \Mp{ \sum_{n=\frac{N}{2}}^{N-100} \n{\widehat{\phi_{n^2}}(\xi)\widehat{K_N^{(3)}}(\xi)}_{L^2} ^q}^{\frac{1}{q}}
    \\ 
    &\leq{}
    C
    \frac{\delta^22^{-N^2}}{N^{\frac{1}{2}}(\log N)^2}
    \Mp{\sum_{n=\frac{N}{2}}^{N-100}2^{qn^2}}^{\frac{1}{q}}\n{F_N}_{L^4}^2\\
    &\leq{}
    C
    \frac{\delta^22^{-200N}}{(\log N)^2}.
\end{align}
Hence, we obtain 
\begin{align}
    \left\{ \sum _{n = \frac{N}{2}}^{N- 100}
    \left( 2^{-\frac{1}{2}n^2} 
    \left\| \psi_{n^2} * \theta_N^{(2)} \right\|_{L^4}
    \right)^q \right\}^{\frac{1}{q}} 
    &\geq{}  
    c
    \left\{ \sum _{n = \frac{N}{2}}^{N- 100}
    \left( 2^{-\frac{1}{2}n^2} 
    \left\| \psi_{n^2} * K_N^{(1)} \right\|_{L^4}
    \right)^q \right\}^{\frac{1}{q}} \\
    &\quad-
    \sum_{m=2}^3
    \left\{ \sum _{n =\frac{N}{2}}^{N- 100}
    \left( 2^{-\frac{1}{2}n^2} 
    \left\| \psi_{n^2} * K_N^{(m)} \right\|_{L^4}
    \right)^q \right\}^{\frac{1}{q}} \\
    &\geq{} c 
    \frac{\delta^2}{(\log N)^2}
    N^{\frac{1}{q}-\frac{1}{2}}
    -
    C
    \frac{\delta^22^{-200N}}{(\log N)^2}\\
    &\geq{} 
    c\frac{\delta^2}{(\log N)^2}
    N^{\frac{1}{q}-\frac{1}{2}}
\end{align}
for sufficiently large $N$.
For the higher iteration part, 
it follows from the same contraction mapping argument as in Step 1 via \eqref{est:1}
that we may construct the solution $\widetilde{\theta}_N$ to \eqref{eq:pQG} with the estimate 
\begin{align}
    \n{\widetilde{\theta}_N}_{\dB_{4,2}^{-\frac{1}{2}}}
    \leq 
    C 
    \frac{\delta^3}{(\log N)^3},
\end{align}
provided that the $\delta$ is sufficiently small.
Hence, the solution $\theta_N = \theta_N^{(1)} + \theta_N^{(2)} + \widetilde{\theta}_N$ satisfies the estimate 
\begin{align}
    \n{\theta_N}_{\dB_{4,q}^{-\frac{1}{2}}}
    \geq{}& 
    \Mp{
    \sum_{n=\frac{N}{2}}^{N-100}
    \sp{
    2^{-\frac{1}{2}n^2}
    \n{\phi_{n^2}*\theta_N}_{L^4}
    }^q
    }^{\frac{1}{q}}\\
    \geq{}&
    c
    \Mp{
    \sum_{n=\frac{N}{2}}^{N-100}
    \sp{
    2^{-\frac{1}{2}n^2}
    \n{\psi_{n^2}*\theta_N^{(2)}}_{L^4}
    }^q
    }^{\frac{1}{q}}
    -
    \n{\theta_N^{(1)}}_{\dB_{4,q}^{-\frac{1}{2}}}
    -
    CN^{\frac{1}{q}-\frac{1}{2}}
    \n{\widetilde{\theta}_N}_{\dB_{4,2}^{-\frac{1}{2}}}\\
    \geq{}&
    c\frac{\delta^2}{(\log N)^2}N^{\frac{1}{q}-\frac{1}{2}}
    -C\frac{\delta}{\log N} 
    -C\frac{\delta^3}{(\log N)^3}N^{\frac{1}{q}-\frac{1}{2}}\\
    ={}&
    \frac{N^{\frac{1}{q}-\frac{1}{2}}}{(\log N)^2}
    \sp{
    c\delta^2
    -C\delta N^{\frac{1}{2}
    -\frac{1}{q}}\log N 
    -\frac{C\delta^3}{\log N}
    }
    \to \infty 
\end{align}
as $N \to \infty$, which completes the proof. 
\end{proof}



\noindent
{\bf Data availability.}\\
Data sharing not applicable to this article as no datasets were generated or analysed during the current study.

\noindent
{\bf Conflict of interest.}\\
The author has declared no conflicts of interest.

\noindent
{\bf Acknowledgements.} \\
The first author was supported by Grant-in-Aid for Research Activity Start-up, Grant Number JP2319011.
The authors would like to express his sincere gratitude to Professor Hideo Kozono for many fruitful advices and continuous encouragement.

\end{document}